\documentclass{svmult}
\usepackage{amsmath,amssymb,amstext}
\usepackage{amsfonts, dsfont}
\usepackage{hyperref}
\usepackage{mathrsfs}
\usepackage[latin1]{inputenc}
\usepackage[T1]{fontenc}

\DeclareMathOperator{\USp}{USp}
\newcommand{\tra}[1]{\,{\vphantom{#1}}^{\mathrm{t}}\hspace{0mm}{#1}}

\newcommand{\1}{\mathds{1}}

\newcommand{\sgn}{\operatorname{{sgn}}}
\DeclareMathOperator{\Prob}{\mathbb{P}}   
\renewcommand{\S}{\mathscr{S}} 
\newcommand{\HH}{{\mathbb H}}                  

\newcommand{\RR}{{\mathbb R}}
\newcommand{\CC}{{\mathbb C}}
 \newcommand{\ii}{{\mathrm{i}}}
  \newcommand{\jj}{{\mathrm{j}}}
  \newcommand{\kk}{{\mathrm{k}}}
 \newcommand{\Id}{{\mathrm{Id}}}
\newcommand{\law}{\overset{\mbox{\rm \scriptsize law}}{=}}

\renewcommand{\Re}{{\mathfrak{Re}}}
\renewcommand{\Im}{{\mathfrak{Im}}}

\title*{Ewens measures  on compact groups and hypergeometric kernels}

\author{\bf P. Bourgade\inst{1}, A. Nikeghbali\inst{2} and A. Rouault\inst{3}}
\institute{Institut Telecom,
46 rue Barrault, 75634 Paris Cedex 13\\ and  Universit\'e Paris 6, LPMA,
175, rue du Chevaleret F-75013 Paris,\\ \textit{e-mail: bourgade@enst.fr}
\and
 Institut f\"ur Mathematik,
 Universit\"at Z\"urich , Winterthurerstrasse 190,
 CH-8057 Z\"urich,
 Switzerland,\\ \textit{e-mail: ashkan.nikeghbali@math.unizh.ch}
\and 
Universit\'e  Versailles-Saint Quentin, LMV,
 B\^atiment Fermat, 45 avenue des Etats-Unis,
78035 Versailles Cedex,\\ \textit{e-mail: alain.rouault@math.uvsq.fr}}

\authorrunning{Bourgade et al.}
\begin{document}

\maketitle
\begin{abstract} 
On 
 unitary compact groups the decomposition of a generic element into product of reflections 
 induces a decomposition of the characteristic polynomial into a product of factors. When the group is equipped with the Haar probability measure, these factors become independent random variables with explicit distributions. Beyond the known results on the orthogonal  and unitary groups ($O(n)$ and $U(n)$), we treat the symplectic case. In $U(n)$, this induces a family of probability changes analogous to the biassing in 
 the Ewens sampling formula known for the symmetric group. Then we study the spectral properties of these measures, connected to the pure Fisher-Hartvig symbol on the unit circle. The associated orthogonal polynomials give rise, as $n$ tends to infinity to a limit kernel at the singularity. 
\end{abstract}
\setcounter{tocdepth}{1}


%
\section{Introduction}
In this paper, $U(n,K)$ is the unitary group over 
 $K=\RR,\CC$ or $\HH$ (the set of real quaternions). 

Let $U$ be distributed with the Haar measure on $U(n,\CC)$. The random variable $\det(\Id_n-U)$ has played a crucial role  in recent years in the study of 
 some  connections between random matrix theory and analytic number theory (see \cite{KeaSna} for more details). In \cite{BHNY}, the authors  show  that $\det(\Id_n-U)$  can be
decomposed as a product of $n$ independent random variables:
\begin{equation}\label{OriginalDecomposition}
\det(\Id_n-U)\law \prod_{k=1}^n
\left(1-\E^{\I\omega_k}\sqrt{B_{1,k-1}}\right),
\end{equation}
where $\omega_1,\dots,\omega_n,B_{1,0},\dots,B_{1,n-1}$ are
independent, 
 the $\omega_k's$ being uniformly
distributed on $(-\pi,\pi)$ and the $B_{1,j}$'s ($0\leq j\leq n-1$)
being beta distributed with parameters 1 and $j$ (with the
convention that $B_{1,0}=1$). In particular, from such a  decomposition, fundamental quantities such as the  Mellin-Fourier transform of $\det(\Id_n-U)$ follow 
 at once. The main ingredient to obtain the decomposition (\ref{OriginalDecomposition}) is 
 a recursive construction of the Haar measure using complex reflections. In particular, every $U\in U(n,\CC)$ can be decomposed as a product of $n$  independent reflections. More precisely, it is proved in \cite{BHNY}  that if  $s_1,\ldots, s_n$ are $n$ independent random variables such that %
for every $k\leq n$, $s_k$ is uniformly distributed on the $k$-th dimensional unit sphere $\S^k$ in $\mathbb C^k$ and if 
$R^{(k)}$ is the reflection of $\CC^k$ mapping $s_k$ onto the first vector of the canonical basis, then
$$R^{(n)} \left(
\begin{array}{cc}
\Id_{1}&0\\
0&R^{(n-1)}
\end{array}
\right) \dots \left(
\begin{array}{cc}
\Id_{n-2}&0\\
0&R^{(2)}
\end{array}
\right) \left(
\begin{array}{cc}
\Id_{n-1}&0\\
0&R^{(1)}
\end{array}
\right) \sim\mu_{U(n,\CC)},$$
where $\mu_{U(n,\CC)}$ stands for the Haar measure on $U(n,\CC)$.
 At this stage two remarks are in order. First,  a similar method works to generate the Haar measure on the orthogonal group $O(n,\RR)$ (see \cite{BHNY}) and this was already noticed by Mezzadri in \cite{mezzadri} using Householder reflections. But as already noticed in \cite{BHNY}, Householder reflections would not work for $U(n,\CC)$ (see next section for more details).  Moreover in \cite{BHNY}, a decomposition such as (\ref{OriginalDecomposition}) could not be 
 obtained for the symplectic group $\USp(2n,\CC)$, which also plays an important role in the connections
between random matrix theory and the study of families of L functions
(see \cite{KatzSarnak1}, \cite{KatzSarnak2}). Indeed, there does not seem to be a natural way to generate recursively the Haar measure on 
this group.\\

\textbf{Question 1.} 
Is there any decomposition  of $\det(\Id_n-U)$ as a product of independent variables of the type (\ref{OriginalDecomposition}), when $U$ is drawn from $\USp(2n,\CC)$, according to the Haar measure?\\

In this paper we shall prove that, in a sense to be made precise, if a subgroup $\mathcal{G}$ of $U(n,K)$ contains enough reflections, then one can recursively generate the Haar measure and obtain a decomposition of the type (\ref{OriginalDecomposition}) for $\det(\Id_n-U)$, $U\in\mathcal{G}$. In particular this will apply to $U(n,\mathbb{H})$ which can be identified with the symplectic group, hence answering question 1 above. Our recursive decomposition of the Haar measure also applies to the symmetric group. This leads us to our second remark concerning the generation of the Haar measure obtained in \cite{BHNY} and explained above. Indeed, this way of generating an element of $U(n,\CC)$ which is Haar distributed by choosing a vector $(s_1,\ldots,s_n)$ of independent variables from $\S^1\times\ldots\times\S^n$, each $s_i$ being uniformly distributed,  is reminiscent of the generation of a random permutation according to the so-called Chinese restaurant process which  we briefly describe  (see \cite{Pit} for a complete treatment). 
Let $[n]$ denote the set $\{1, \cdots, n\}$ and ${\mathcal S}_n$ the symmetric group of order $n$.
It is known that for $n \geq 2$, every permutation $\sigma \in {\mathcal S}_n$ can be decomposed in the following way:
\begin{equation}
\label{dec0}
\sigma = \tau_n \circ \dots \circ \tau_2
\end{equation}
where for $k = 2, \dots, n$,  either $\tau_k$ is the identity or  $\tau_k$ is the transposition $(k, m_k)$ for some $m_k \in [k-1]$. In the first case we will say by extension that it is the transposition $(k,m_k)$ with $m_k = k$. 
This decomposition is unique, see Tsilevich \cite{Tsilevich}, the lemma p. 4075. It corresponds to the Chinese restaurant generation of a permutation. Let us consider cycles as "tables". 
Integer $1$ goes to the first table. If $\tau_2 \not= \hbox{Id}$, then integer $2$ goes to the first table, at the left of $1$. If $\tau_2= \hbox{Id}$, it goes to a new table. When integers $1, \dots, k$ are placed, then $k+1$ goes to a new table if $\tau_{k+1}=\hbox{Id}$, and goes to the left of $\tau_{k+1}(k+1) =j_{k+1}$ if not. We get a bijection between $[1]\times [2]\times \cdots \times [n] \rightarrow {\mathcal S}_n$.
It is projective (or consistent) in the sense that if $\sigma$ is in ${\mathcal S}_{n+1}$ the restriction of $\sigma$ to $[n]$ is in ${\mathcal S}_n$.

In this setting, the number of cycles $k_\sigma$ of a permutation $\sigma$ is the number of tables, i.e. the number of $\Id$ in (\ref{dec0}) i.e. \footnote[1]{The other construction of a random permutation named Feller's coupling uses the variables in the reverse order $\xi_n, \cdots, \xi_1$, but this construction is not projective.}
\begin{equation}\label{ksigma}k_\sigma = \sum_1^n \xi_r\,,\end{equation}
where $\xi_r = 1(\tau_r = \hbox{Id})$. For a matricial rewriting, we make a change of basis. Let $e'_j = e_{n-j+1}$ and let $R^{(k)}$ be the restriction of $\tau_k$ to $[k]$. Then the product in (\ref{dec0}) is represented by
$$
R^{(n)} \left(
\begin{array}{cc}
\Id_1&0\\
0&R^{(n-1)}
\end{array}
\right) \dots \left(
\begin{array}{cc}
\Id_{n-2}&0\\
0&R^{(2)}
\end{array}
\right).$$
 If at each stage, the integer $m_k$ is chosen uniformly in $[k]$, then the induced measure on ${\mathcal S}_n$ is the uniform distribution denoted by $\mu_{\mathcal{S}_n}$.

Actually, one can more generally generate in this way the Ewens measure on 
 $\mathcal{S}_n$ (see Tsilevich \cite{Tsilevich} and Pitman \cite{Pit}). 
The Ewens measure $\mu^{(\theta)}$, $\theta>0$, is a deformation of 
 $\mu_{\mathcal{S}_n}$ 
  obtained by performing a change of probability measure or a sampling 
   in the following way:
\begin{equation}\label{ew11}
\mu^{\theta}_n(\sigma)=\dfrac{\theta^{k_\sigma}}{(\theta)_n} \cdot \mu_{\mathcal{S}_n}(\sigma)\,.
\end{equation}
To generate $\mu^{\theta}_n$, 
 one has to pick $n$ integers $m_1,m_2,\ldots,m_n$, independently,  from 
$[1]\times \cdots\times[n]$ 
  according to the probability distribution 

\begin{eqnarray*}
\mathbb{P}(m_k=k)= \dfrac{\theta}{\theta+k-1} \ , \ 
\mathbb{P}(m_k=j) = \dfrac{1}{\theta+k-1} \ j = 1, \cdots,  k-1\,.
\end{eqnarray*}

\textbf{Question 2.} Is there an analogue of the Ewens measure on the unitary group $U(n, \CC)$?\\

We shall see  in this paper that there indeed exists an analogue of the Ewens measure on  $U(n, \CC)$: more precisely we generalize (\ref{ew11}) to unitary groups and a
particular class of their subgroups. The analogue of transpositions
 are reflections and the weight of the sampling  is now 
$\det(\Id-U)^{\overline{\delta}}\det(\Id-\overline{U})^{\delta}$,
$\delta\in\CC$, $\Re(\delta)>-1/2$,
so that the measure
$\mu^{(\delta)}_{U(n)}$  on $U(n)$, which is defined by
$$
\mathbb{E}_{\mu^{(\delta)}_{U(n)}}\left(f(U)\right)=\frac{\mathbb{E}_{\mu_{U(n)}}\left(f(U)\det(\Id-U)^{\overline{\delta}}
\det(\Id-\overline{U})^{\delta}\right)}{\mathbb{E}_{\mu_{U(n)}}\left(\det(\Id-U)^{\overline{\delta}}
\det(\Id-\overline{U})^{\delta}\right)}
$$
for any  
test function $f$, is the analogue of the Ewens measure. 
Such samplings with $\delta\in\RR$ have already been studied on the
finite-dimensional unitary group by Hua \cite{Hua}, and results
about the infinite dimensional case (on complex Grassmannians) were
given by Pickrell (\cite{Pick1} and \cite{Pick2}). More recently,
Neretin \cite{Ner} also considered this 
 measure, introducing
the possibility $\delta\in\CC$. Borodin and Olshanski \cite{BO}
have used the analogue of this 
 measure in the framework of  the infinite
dimensional unitary group and proved ergodic properties.
Forrester and Witte in \cite{Forr} referred to this measure
as the cJUE distribution. 
We also studied this ensemble in \cite{BNR} in relation with the theory of orthogonal polynomials on the unit circle.  Following \cite{Forr} and \cite{BNR} we shall call the ensemble of unitary matrices endowed with this sampled measure the \emph{circular Jacobi ensemble}. 

It is natural to ask whether the circular Jacobi ensemble has some interesting properties: indeed, the case $\delta=0$ corresponds to the Haar measure and it is well known this ensemble enjoys many remarkable spectral properties. For instance, the point process associated to the eigenvalues  is determinantal and the associated rescaled kernel converges to the 
 sine kernel .  
The projection of
the  measures $\mu^{(\delta)}_{U(n)}$ on the spectrum has the density
\[\frac{1}{{\mathcal Z}_n} \prod_{j=1}^n w^{\mathbb T}(\E^{\I\theta_j}) \prod_{1\leq i<j\leq n}|\E^{\I\theta_i} - \E^{\I\theta_j}|^2\]
where the weight $w^{\mathbb T}$  on $\mathbb T  = \{ \E^{\I\theta} , \theta \in [-\pi, \pi]\}$ is defined by
\[w^{\mathbb T} (\E^{\I\theta}) = (1 - \E^{\I\theta})^{\bar \delta}(1 - \E^{-\I\theta})^{\delta} = (2-2\cos\theta)^a \E^{-b(\pi \sgn\theta-\theta)}\,,\]
$(\delta = a + \I b)$ and ${\mathcal Z}_n$ is a normalization constant. 
Note that when $b \not = 0$,
 an asymmetric singularity at $1$ occurs. The statistical
properties of the $\theta_k$'s depend on the successive orthonormal
polynomials $(\varphi_k)$ with respect to the normalized version $\widetilde w^{\mathbb T}$ of $w^{\mathbb T}$ and 
%
the normalized reproducing kernel 
\[\widetilde K^{\mathbb T}_n(\E^{\I \theta}, \E^{\I\tau})= \sqrt{\widetilde w^{\mathbb T} (\E^{\I \theta})\widetilde w^{\mathbb T} (\E^{\I\tau})} \sum_{\ell =0}^{n-1} \overline{\varphi_\ell(\E^{\I \theta})}\varphi_\ell (\E^{\I\tau}) 
\,.
\]

In \cite{BO} the authors consider the image of $\mu_{U(n)}^{(\delta)}$ 
 by the Cayley transform on the set of Hermitian matrices and make a thorough study of the spectral properties of this random matrix ensemble. In particular they prove that the eigenvalues  form a determinantal process and show that the associated rescaled kernel converges to some hypergeometric kernel. As expected, we shall see that the eigenvalues process of the circular Jacobi ensemble is also determinantal and for every $n$, we identify 
 the  hypergeometric kernel $K_n^{(\delta)}$ associated with it.  \\

\textbf{Question 3.} Is there an appropriate rescaling of  the kernels $K_n^{(\delta)}$ such that the rescaled kernels converge to some kernel $K_\infty^{(\delta)}?$ \\

We shall see that the answer to question 3 is positive and that  the kernel $K_\infty^{(\delta)}$ is  a \textsl{confluent hypergeometric kernel}, 
with a natural connection to that obtained by Borodin and Olshanski in \cite{BO} on the set of Hermitian matrices.
The case $\delta=0$ corresponds to the sine kernel. 

The weight $w^{\mathbb T}$ is a generic example leading to a
singularity
$$
c^{(+)}|\theta|^{2a}\1_{\theta>0}+c^{(-)}|\theta|^{2a}\1_{\theta<0}
$$
at $\theta=0$, with distinct positive constants $c^{(+)}$ and
$c^{(-)}$. The confluent hypergeometric kernel, depending on the two
parameters $a$ and $b=\frac{1}{2\pi}\log(c^{(-)}/c^{(+)})$, is
actually universal for the measures presenting  the above
singularity, as proved in a forthcoming paper, following the method initiated by Lubinsky (\cite{Lubinsky1}, \cite{Lubinsky2}).  For a universality result when $\delta$ is real see \cite{RS}.

The layout of the paper is as follows. In Section \ref{Gen} we present the generation by reflections and  deduce a splitting formula for the characteristic polynomial (Theorem \ref{thm:CompleteDecDet}). As an application, we define the generalized Ewens measure depending on the complex parameter $\delta$ (Theorem \ref{thm:GeneralEwens}).   Section \ref{Hyper} is devoted to a study of the kernel which governs the correlations of eigenvalues when the unitary group is equipped with this measure and its asymptotics (Theorem \ref{theokernel}). The main properties of the families of hypergeometric functions $_2F_1$ and $_1F_1$ are recalled in the Appendix. 

\section{Generating the Haar measure and the generalized Ewens measure}
\label{Gen}
\subsection{Complex reflections}
Reflections play a central role in the generation of the Haar measure for the classical compact groups. In the case of $O(n)$ the decomposition into a product of reflections is well known, see \cite{Diaconis} and other references as explained in \cite{mezzadri}. Householder reflections are generally used  in the case of $O(n)$, but they are not suitable for $U(n,\CC)$. 
Indeed, recall that Householder reflections are of the form $H_v= \Id -2 v\langle v|\ \cdot\rangle$. For every unit $y$, it is possible to choose $v$ such that 
$H_v y = \alpha e_1$ with $\alpha = \pm\frac{y_1}{|y_1|}$, where $e_1$ is the first element of the canonical basis. So when the ground field is $\CC$, then $\alpha\neq1$ in general and there does not exist a Householder reflection which maps $y$ onto $e_1$, whereas this can always be achieved when the ground field is $\RR$. That is why it is not possible to directly extend the arguments in \cite{mezzadri} to $U(n,\CC)$. In \cite{BHNY} and \cite{BNR} it is proposed to use complex (resp. quaternionic) proper reflections, that is  norm preserving automorphisms of $\mathbb C^n$ (resp. $\mathbb H^n$) that leave exactly one hyperplane pointwise fixed. So 
a reflection will be either the identity or a unitary transformation $U$ such that $I-U$ is of rank one. It may be written as
\[s_{a, \lambda} (y) = y - a\frac{(1-\lambda) \langle a, y\rangle}{|a|^2}\]
where $a \in \mathbb H^n$ and $\lambda \in \mathbb H$ with $|\lambda|=1$ ($\lambda$ is the second eigenvalue). 
If $x \not=e_1$, there exists a  reflection mapping $e_1$ onto $x$. It is enough to take
$a = e_1 - x$ and $\lambda = -(1-x_1)(1-\bar x_1)^{-1}$ where $x_1 = \langle e_1, x\rangle$.

\subsection{Generating the Haar measure on $U(n,K)$ and on some of its subgroups}

We first give conditions under which an element of a subgroup of  $U(n,K)$
(under the Haar measure) can be generated as
a product of independent reflections. This will lead
to some remarkable identities for the characteristic polynomial.

Let
$(e_1,\dots,e_n)$ be an orthonormal basis of $\mathbb K^n$. Let $\mathcal{G}$ be a subgroup
of $U(n,K)$ and for
all $1\leq k\leq n-1$, let $$\mathcal{H}_k=\{G\in\mathcal{G}\mid G(e_j)=e_j, \ 1 \leq j \leq k\}\,,$$ the subgroup of
$\mathcal{G}$ which stabilizes $e_1, \cdots, e_k$. We set ${\mathcal H}_0 = {\mathcal G}$. For a generic compact group
$\mathcal{A}$, we write $\mu_{\mathcal{A}}$ for the unique Haar
probability measure on $\mathcal{A}$. Finally for
all $1\leq k\leq n$ let $p_k$ be the map $U\mapsto
U(e_k)$.

\begin{proposition}\label{prop:decompositionRank1}
Let $G\in\mathcal{G}$ and $H\in\mathcal{H}_1$ be independent random matrices, and assume that 
 $H \sim \mu_{\mathcal{H}_1}$. Then $GH
\sim\mu_\mathcal{G}$ if and only if $G(e_1)\sim
p_1(\mu_{\mathcal{G}})$.
\end{proposition}
\begin{proof}
The proof is exactly the same as in \cite{BHNY} Prop. 
 2.1, changing $U(n+1)$ into $\mathcal G$ and $U(n)$ into $\mathcal H$.
\end{proof}

\begin{definition}\label{coherent}
A sequence $(\nu_0,\dots,\nu_{n-1})$ of probability measures on
$\mathcal{G}$ is said to be coherent with $\mu_\mathcal{G}$ if for
all $0\leq k\leq n-1$, 
\[\nu_k(\mathcal{H}_{k})=1\ \hbox{and}\ p_{k+1}(\nu_{k})= p_{k+1}(\mu_{{\mathcal H}_k})\,.\]
\end{definition}
In the following, $\nu_0\star\nu_1\star\dots\star\nu_{n-1}$
stands for the law of a random variable $H_0H_1\dots H_{n-1}$ where
all $H_i$'s are independent and $H_i\sim \nu_i$. Now we can provide
a general method to generate an element of $\mathcal{G}$ endowed
with its Haar measure.
\begin{theorem}\label{thm:decompositionRankN}
If $\mathcal{G}$ is a subgroup of $U(n, K)$ and 
$(\nu_0,\dots,\nu_{n-1})$ is a sequence of coherent measures with
$\mu_\mathcal{G}$, then we have:
$$\mu_\mathcal{G}=\nu_0\star\nu_1\star\dots\star\nu_{n-1}.$$
\end{theorem}

\begin{proof}
It is sufficient to prove by induction on $1\leq k\leq n$ that
$$
\nu_{n-k}\star\nu_{n-k+1}\star\dots
\star\nu_{n-1}=\mu_{\mathcal{H}_{n-k}},
$$
which gives the desired result for $k=n$. If $k=1$ this is obvious.
If the result is true at rank $k$, it remains true at rank $k+1$ by
a direct application of Proposition \ref{prop:decompositionRank1} to
the groups $\mathcal{H}_{n-k-1}$ and its subgroup
$\mathcal{H}_{n-k}$.
\end{proof}

As an example, take the orthogonal group $O(n)$. Let
$\S^{(k)}_\RR$ be the unit sphere $\{x\in\RR^k\mid |x|=1\}$ and, for
$s_k\in\S^{(k)}_\RR$, let $R^{(k)}$ be the matrix of the
reflection which transforms $s_k$ into $e_1$.
If  $s_k$ is uniformly distributed on  $\S^{(k)}_\RR$ and if all the $s_k$ are independent,
 then by Theorem \ref{thm:decompositionRankN}, the matrix 
$$
R^{(n)} \left(
\begin{array}{cc}
1&0\\
0&R^{(n-1)}
\end{array}
\right) \dots \left(
\begin{array}{cc}
\Id_{n-2}&0\\
0&R^{(2)}
\end{array}
\right) \left(
\begin{array}{cc}
\Id_{n-1}&0\\
0&R^{(1)}
\end{array}
\right).$$
is $\mu_{O(n)}$ distributed.

\subsection{Splitting of the characteristic polynomial}
In view to phrase a general version of formula (\ref{OriginalDecomposition}) which is proved in \cite{BHNY}, we need the following definition:
\begin{definition}\label{defn:DensityOfReflections}
Note $\mathcal{R}_k$ the set of elements in $\mathcal{H}_k$ which
are reflections. If for all $0\leq k\leq n-1$
$$
\{R(e_{k+1})\mid R\in\mathcal{R}_k\}=\{H(e_{k+1})\mid
H\in\mathcal{H}_k\},
$$
the group $\mathcal{G}$ will be said to satisfy condition (R) (R
standing for reflection).
\end{definition}
\begin{remark}
It is easy to see  that $U(n,K)$ and $\mathcal{S}_n$  satisfy condition (R). In the next subsection we shall see more examples.
\end{remark}
\begin{lemma}\label{lem1}
Let $\mathcal{G}$ be a subgroup of $U(n,K)$ which satisfies condition (R). Let $G\in\mathcal{G}$. Then  there exist 
 reflections $R_k\in\mathcal{R}_k$, $0\leq k\leq n-1$, such that
\begin{equation}\label{decoref}
G=R_0R_1\ldots R_{n-1}.
\end{equation}  
\end{lemma}
\begin{proof}
This result has been established in \cite{BNR} when $\mathcal{G}=U(n,\CC)$. The proof in this more general case goes exactly along the same line.
\end{proof}
The following deterministic lemma is a key result to obtain a decomposition of $\det(\Id_n -U)$ as a product of independent random variables:
\begin{lemma}
\label{determ}
If for $k=1, \dots, n-1$, $R_k \in {\mathcal R}_k$, then
\begin{equation}
\det (\Id_n - R_0\cdots R_{n-1}) = \prod_{k=0}^{n-1} \left(1 - \langle e_{k+1}, R_k(e_{k+1}\rangle\right)\,.
\end{equation}
\end{lemma}
\proof
We start with $\det(\Id_n -RH) = (\det H) \det (H^* - R)$.
Since $H$ (hence $H^*$), stabilizes $e_1$, we have 

i) $(H^* -R)(e_1) = e_1 - R(e_1) =: a$ (say), 

ii) for $w\perp e_1$, $H^* (w) \perp e_1$ and since $R$ is a reflection, $R(w) - w$ is a scalar multiple of $a$. 

By the multilinearity of the determinant, we get
\[
\det(H^* - R) = \langle e_1 , e_1 - R(e_1)\rangle \det (\pi(H^*) -\Id_{n-1})\,
\]  
which yields
\[\det(\Id_n - RH) = (1- \langle e_1, R(e_1)\rangle) \det(\Id_{n-1} -\pi(H))\,.\]
Iterating, we can conclude.
\qed

The following result now follows immediately from Theorem \ref{thm:decompositionRankN} and Lemmas \ref{lem1} and \ref{determ}.

\begin{theorem}\label{thm:CompleteDecDet}
Let $\mathcal{G}$ be a subgroup of $U(n,K)$ satisfying condition (R), and let
$(\nu_0,\dots,\nu_{n-1})$ be coherent with $\mu_{\mathcal{G}}$. 
If $G\sim\mu_\mathcal{G}$, then
$$
\det(\Id-G)\law\prod_{k=0}^{n-1}\left(1-\langle
e_{k+1}, H_k(e_{k+1}),\rangle\right).
$$
where $H_k\sim\nu_k$, $0\leq k\leq n-1$, are independent.  
\end{theorem}

\subsection{Applications}
\subsubsection{The symmetric group.}

Consider now $\mathcal{S}_n$ the group of permutations of size $n$. An
element $\sigma\in\mathcal{S}_n$ can be identified with the matrix
$(\delta^j_{\sigma(i)})_{1\leq i,j\leq n}$ ($\delta$ is Kronecker's
symbol). It is clear that $1$ is eigenvalue of this matrix, with eigenvector $e_1+ \cdots + e_n$. Ben Hambly et al. \cite{BenHambly} considered the characteristic polynomial at $s\not=1$. To make relevant our problem of determinant splitting, we introduce wreath products, following the definition of Wieand  \cite{wieand}.  

Let $F$ be a subgroup of $\mathbb T = \left\{x\in \mathbb C\mid |x|^2=1 \right\}$, endowed
with the Haar probability measure $\mu_F$.
 Then
the wreath product 
$F\wr \mathcal{S}_n$ provides 
another example of determinant-splitting. 
 An element of $F^n$  can be thought of as a function
from the set $[n]$ to $F$. The group ${\mathcal S}_n$  acts on $F^n$ in the
following way: if $f=(f(1),\dots, f(n)) \in  F^n$ and $\sigma \in  {\mathcal S}_n$,  define $
f_\sigma\in F^n$  to be the function $f_\sigma =f\circ \sigma^{-1}$. Finally take the product on $F^n$ to
be $(f(1),\dots, f(n))\cdot(g(1),\dots, g(n))=(fg(1),\dots, fg(n))$. The wreath product of $F$ by $\mathcal{S}_n$ , denoted $F\wr \mathcal{S}_n$ , is the
group of elements $\{(f; s): f \in F^n, \sigma\in \mathcal{S}_n\}$ with multiplication 
\[(f; \sigma)\cdot(h; \sigma')
=(fh_\sigma; \sigma\sigma')\,.\]
If we  represent  $(f; \sigma)$ by the matrix $(f(i)\delta^i_{\sigma(j)})_{1\leq i,j\leq n}$, then the product in $F\wr \mathcal{S}_n$ corresponds to the usual matricial product which makes $F\wr \mathcal{S}_n$ a subgroup of $U(n, \mathbb C)$. The usual examples are $F= \{1\}$, $F= \mathbb Z_2$ and $F = \mathbb T$.

\begin{corollary} \label{permutationscase}
Let $G$ $\in \mathcal{G}(= F\wr \mathcal{S}_n)$ be
$\mu_{\mathcal{G}}$ distributed. Then
$$
\det (\Id_n-G)\law \prod_{j=1}^n \left(1-\varepsilon_j X_j\right),
$$
with $\varepsilon_1,\dots,\varepsilon_n,X_1,\dots,X_n$ independent random variables, the
$\varepsilon_j$'s $\mu_F$ distributed, $\mathbb P(X_j=j)=1/j$,
$\mathbb P(X_j=0)=1-1/j$.
\end{corollary}

\proof We apply Theorem \ref{thm:CompleteDecDet}.
As  reflections  correspond now to  transpositions, condition
(R) holds. Moreover $R_k(e_{k+1})$ is uniformly distributed on
the  set 
 $Fe_{k+1} \cup \dots\cup Fe_n$, so that $\langle e_{k+1}, R_k(e_{k+1})\rangle$ is $0$ with probability $(n-k)/n$ and otherwise, it is uniform on $F$.
 \qed

\begin{remark} 
Notice that if $G = (f;\sigma)$ with $\sigma = \tau_n\circ \cdots\circ\tau_2$ (cf. (\ref{dec0})), then $X_j$ is the indicator function of $\tau_{n-j+1} = \Id$.
\end{remark}

\subsubsection{Unitary and orthogonal groups}
Take $\mathcal{G}=U(n,\CC)$. Then $\mu_{\mathcal{H}_k}
=f_k(\mu_{U(n-k,\CC)})$ where $f_k : A\in U(n-k,\CC)\mapsto \Id_{k}
\oplus A$. As all reflections with respect to a hyperplane of
$\CC^{n-k}$ are elements of $U(n-k,\CC)$, one can apply Theorem \ref{thm:decompositionRankN} and Lemma \ref{determ}. The Hermitian products $\langle
e_k,h_k(e_k)\rangle$ are distributed as the first coordinate of the
first vector of an element of $U(n-k,\CC)$, that is to say the first
coordinate of the $(n-k)$-dimensional unit complex sphere with
uniform measure : \[\langle e_{k+1}, H_k(e_{k+1})\rangle\law
\E^{\I\omega_n}\sqrt{B_{1,n-k-1}}\] with $\omega_n$ uniform on
$(-\pi,\pi)$ and independent of $B_{1,n-k-1}$, a beta variable with
parameters 1 and $n-k-1$.

Therefore, as a consequence of Theorem
\ref{thm:CompleteDecDet}, we obtain the following decomposition formula derived in  \cite{BHNY}.
 For 
  $g \in U(n,\CC)$ which is  $\mu_{U(n,\CC)}$
distributed, one has 
$$
\det (\Id_n-G)\law \prod_{k=1}^n
\left(1-\E^{\I\omega_k}\sqrt{B_{1,k-1}}\right),
$$
with $\omega_1,\dots,\omega_n,B_{1,0},\dots,B_{1,n-1}$ independent
random variables, the $\omega_k$'s uniformly distributed on
$(-\pi,\pi)$ and the $B_{1,j}$'s ($0\leq j\leq n-1$) being beta
distributed with parameters 1 and $j$ (by
convention, $B_{1,0}=1$).\\

 A similar reasoning may be applied to $SO(2n)$ 
(with the complex unit spheres replaced by the real ones) to yield the following:
let $G\in SO(2n)$ be $\mu_{SO(2n)}$ distributed, then (Corollary 6.2 in \cite{BHNY})
$$
\det (\Id_{2n}-G)\law 2 \prod_{k=2}^{2n}
\left(1-\epsilon_k\sqrt{B_{\frac{1}{2},\frac{k-1}{2}}}\right)\,.
$$
\subsubsection{The quaternionic group}
Our goal with this example is to solve Question 1 which was raised in the Introduction. To this end we establish an analogous to 
 Lemma \ref{determ} 
 and 
  use the fact that  $U(n,\mathbb H) \cong \USp(2n)$ which is also denoted $Sp(n)$,  see for instance \cite{mezzadri} Theorem 2. Then we apply Theorem \ref{thm:decompositionRankN}. Let us give details. Recall that the symplectic group $\USp(2n,\CC)$ is defined as  $\USp(2n,\CC)=\{U\in U(2n,\CC)\mid
UJ_n\tra{U}=J_n\}$, with
\begin{equation}\label{J}
J_n=\left(\begin{array}{cc}0&\Id_n\\-\Id_n&0\end{array}\right).
\end{equation}
Let $$
\phi : \left\{
\begin{array}{ccc}
\HH&\to&M(2,\CC)\\
a+\ii b+\jj c+\kk d&\mapsto& \left(
\begin{array}{cc}
a+\ii b&c+\ii d\\
-c+\ii d&a-\ii b
\end{array}
\right)
\end{array}
\right.,
$$
be the usual representation of quaternions. It is a continuous injective ring morphism 
such that $\phi(\bar x) = \phi(x)^*$. It induces the ring morphism
$$
\Phi : \left\{
\begin{array}{ccc}
M(n,\mathbb H)&\to&M(2n,\mathbb C)\\
(a_{ij})_{1\leq i,j\leq n}&\mapsto& (\phi(a_{ij}))_{1\leq i,j\leq n}
\end{array}
\right..
$$
In particular
\[
\Phi(U(n,\mathbb H)) = \{ G \in U(2n,\CC) : G\tilde{Z}_n\tra{G}=\tilde{Z}_n\}\] 
where
$\tilde{Z}_n=J_1\oplus\dots\oplus J_1$ and  
$J_1=\left(
\begin{array}{cc}
0\ &1\\
-1\ &0
\end{array}
\right)$.  Since $\tilde{Z}_n$ is conjugate to
$J_n$, defined by (\ref{J}), the set $\Phi(U(n,\mathbb H))$ is therefore
conjugate to $\USp(2n,\CC)$.
We can therefore consider $\det (I- \Phi(G))$

\begin{lemma}
\label{determsympl}
If for $k=1, \dots, n-1$, $R_k \in {\mathcal R}_k$, then
\begin{equation}
\det (\Id_{2n} - \Phi(R_0\cdots R_{n-1})) = \prod_{k=0}^{n-1} \det(\Id_2 - \phi(\langle e_{k+1}, R_k(e_{k+1}\rangle))\,.
\end{equation}
\end{lemma}

\proof
Let us first remark that the canonical basis $e_1, \dots, e_n$ of $\mathbb H^n$ is mapped by $\Phi$ into the canonical basis $\varepsilon_1, \dots, \varepsilon_{2n}$ of $\mathbb C^{2n}$, where the $2n\times 2$ matrix $[\varepsilon_{2k-1},\varepsilon_{2k}]$ is exactly  $\Phi(e_k)$. Moreover, if $R$ is a proper reflection (leaving invariant an hyperplane), $\Phi(R)$ is a bireflection of $\mathbb C^{2n}$ i.e. a unitary transformation leaving invariant a vector space of codimension $2$.

We start with \[\det\big((\Id_{2n} -\Phi(RH)\big ) = \det\big(\Id_{2n} -\Phi(R)\Phi(H)\big)= \det \Phi(H) \det \big(\Phi(H^*) - \Phi(R)\big)\]
Since $H$ (hence $H^*$)  stabilizes $e_1$, then $\Phi(H)$ (and $\Phi(H)^*$) stabilizes $\varepsilon_1$ and $\varepsilon_2$, so we have:

i) $(H^* -R)(e_1) = e_1 - R(e_1) =: a = [a_1, a_2]$ (say), 
hence, for $i=1,2$, $$(\Phi(H^*) -\Phi(R))(\varepsilon_i) = \varepsilon_i - \Phi(R)(\varepsilon_i) =: a_i\,.$$ 

ii) Assume that $\langle e_1, w\rangle = 0$.  Trivially, $\langle e_1, H^* (w)\rangle =0$ hence
$\Phi(H^*)(w)$ is a matrix whose column vectors are orthogonal to $\varepsilon_1$ and $\varepsilon_2$. 
Moreover, since $R$ is a quaternionic reflection, $R(w) - w$ is a (right) scalar multiple of $a$ (see \cite{CohenAM} Proposition 1.6),  so
$\Phi\left(R(w) - w\right)$ is a $2n\times 2$ matrix whose columns are in Span $(a_1, a_2)$.

By the multilinearity of the determinant, we get
\[
\det\big(\Phi(H)^* - \Phi(R)\big) = \det\left(\langle \epsilon_i, a_j\rangle_{1\leq i,j\leq 2}\right)
 \det (\pi(H^*) -\Id_{2n-2})\,
\]  
which yields
\[\det(\Id_n - \Phi(RH)) = \det(\Id_2- \phi(\langle e_1, R(e_1)\rangle)) \det(\Id_{2n-2} -\pi(H))\,.\]
Iterating, we can conclude. \qed

\begin{corollary}\label{symplectic}{\it\bf Symplectic group.}
Let $G\in \USp(2n,\CC)$ be $\mu_{\USp(2n,\CC)}$ distributed. Then
$$
\det (\Id_{2n}-G)\law \prod_{k=1}^{n}
\left((a_k-1)^2+b_k^2+c_k^2+d_k^2\right),
$$
where the vectors $(a_k,b_k,c_k,d_k)$, $1\leq k\leq n$ are 
independent and $(a_k,b_k,c_k,d_k)$ are 4 coordinates of the
4k-dimensional real unit sphere endowed with the uniform measure.
\end{corollary}

\begin{remark}
We have $(a_k,b_k,c_k,d_k)\law\frac{1}{\sqrt{\mathcal{N}_1^2+\dots+\mathcal{N}_{4k}^2}}(\mathcal{N}_1,\mathcal{N}_2,\mathcal{N}_3,\mathcal{N}_4)$,
with the $\mathcal{N}_i's$ i.i.d. ${\mathcal N}(0,1)$.
Now, since for $p< q$
$$\frac{\mathcal{N}_1^2+\dots+\mathcal{N}_p^2}{\mathcal{N}_1^2+\dots+\mathcal{N}_q^2}\law B_{\frac{p}{2},\frac{q-p}{2}}\,,$$
\newpage
 we get 
 the somehow more tractable identity in law
$$
\det (\Id_{2n}-G)\law \prod_{k=1}^{n}
\left(\left(1+\epsilon_k\sqrt{B_{\frac{1}{2},2k-\frac{1}{2}}}\right)^2+\left(1-B_{\frac{1}{2},2k-\frac{1}{2}}\right)B'_{\frac{3}{2},2k-2}\right),
$$
with all variables independent,
$\Prob(\epsilon_k=1)=\Prob(\epsilon_k=-1)=1/2$.

This method can be applied to other interesting
groups such as $\USp(2n,\RR)=\{u\in U(2n,\RR)\mid uz\tra{u}=z\}$
thanks to the morphism
$$
\phi : \left\{
\begin{array}{ccc}
\CC&\to&M(2,\RR)\\
a+\ii b&\mapsto& \left(
\begin{array}{cc}
a&-b\\
b&a
\end{array}
\right)
\end{array}
\right..
$$
The traditional representation of the quaternions in $M(4,\RR)$
$$
\phi : \left\{
\begin{array}{ccc}
\CC&\to&M(4,\RR)\\
a+\ii b+\jj c+\kk d&\mapsto& \left(
\begin{array}{cccc}
a&-b&-c&-d\\
b&a&-d&-c\\
c&d&a&-b\\
d&-c&b&a
\end{array}
\right)
\end{array}
\right.
$$
gives another identity in law for a compact subgroup of
$U(4n,\RR)$.
\end{remark}

\subsection{The generalized Ewens measure}
In this section we wish to define a generalization of the Ewens measure on $U(n,K)$ and some of its subgroups which will agree with the classical definition on the symmetric group. We first recall the definition of the Ewens measure on the symmetric group and how it can be generated. 
\subsubsection{The Ewens measure on $\mathcal{S}_n$}
Recall (see (\ref{dec0}) Section 1) 
 that every permutation $\sigma \in {\mathcal S}_n$ can be decomposed in the following way:
\begin{equation}
\label{dec1}
\sigma = \tau_n \circ \dots \circ \tau_2
\end{equation}
where for $k = 2, \dots, n$,   $\tau_k$ is either the identity or  
 the transposition $(k, m_k)$ for some $m_k \in [k-1]$. In the first case we will say by extension that it is the transposition $(k,m_k)$ with $m_k = k$.  
 The number of cycles in the decomposition of $\sigma$ is denoted $k_\sigma$.
The system of Ewens measures of parameter $\theta>0$ consists in choosing the $m_k, k=1, \dots, n$ independently, with distribution
$$\mathbb P(m_k =k) = \frac{\theta}{\theta +k -1} \ ; \ \mathbb P(m_k = j) = \frac{1}{\theta +k-1} , j =1, \dots, k-1.$$
It is known that the induced probability on ${\mathcal S}_n$ is
\begin{equation}\label{ew34}
\mu_{n}^{\theta}(\sigma)= \frac{\theta^{k_\sigma}}{(\theta)_n}.
\end{equation}
\subsubsection{The generalized Ewens measure}
In the following, $\mathcal{G}$ is any subgroup of $U(n,K)$. Take $\delta\in\CC$ such that
\begin{equation}\label{finiteExpectation}
0<\mathbb E_{\mu_\mathcal{G}}\left(\det(\Id_n-G)^{\overline{\delta}}\det(\Id_n-\overline{G})^\delta\right)<\infty.
\end{equation}
For $0\leq k\leq n-1$ we note
$$
\exp^{(k)}_\delta:\left\{
\begin{array}{ccl}
\mathcal{G}&\to& \RR^+\\
G&\mapsto& (1-\langle
e_{k+1},G(e_{k+1})\rangle)^{\overline{\delta}}(1-\langle
\overline{e_{k+1},G(e_{k+1})\rangle})^{\delta}
\end{array}
\right..
$$
Moreover, define $\det_\delta$ as the function
$$
\mbox{det}_\delta:\left\{
\begin{array}{ccl}
\mathcal{G}&\to& \RR^+\\
G&\mapsto&
\det(\Id_n-G)^{\overline{\delta}}\det(\Id_n-\overline{G})^\delta
\end{array}
\right..
$$
 Then the following generalization of Theorem \ref{thm:decompositionRankN}
 (which corresponds to the case $\delta=0$) holds.
 However, note that, contrary to Theorem \ref{thm:decompositionRankN},
 in the following result we need that the coherent measures 
  be supported
 by the set of reflections.

\begin{theorem}{\rm \bf Generalized Ewens sampling formula.} \index{Ewens sampling formula! on general compact groups} \label{thm:GeneralEwens}
Let $\mathcal{G}$ be a subgroup of $U(n, K)$ checking condition (R)
and (\ref{finiteExpectation}). Let $(\nu_0,\dots,\nu_{n-1})$ be a
sequence of measures coherent with $\mu_\mathcal{G}$, with
$\nu_k(\mathcal{R}_k)=1$. We note $\mu^{(\delta)}_\mathcal{G}$ the
$\det_\delta$-sampling of $\mu_\mathcal{G}$ and $\nu^{(\delta)}_k$
the $\exp^{(k)}_\delta$-sampling of $\nu_k$. Then
$$\nu^{(\delta)}_0\star\nu^{(\delta)}_1\star\dots\star\nu^{(\delta)}_{n-1}=\mu^{(\delta)}_\mathcal{G},$$
i.e., for all test functions $f$ on
$\mathcal{G}$,
$$
\mathbb E_{\nu^{(\delta)}_0\star\dots\star\nu^{(\delta)}_{n-1}}\left(f(R_0R_1\dots
R_{n-1})\right)=
\frac{\mathbb E_{\mu_{\mathcal{G}}}\left(f(G)\det(\Id_n-G)^{\overline{\delta}}\det(\Id_n-\overline{G})^\delta\right)}
{\mathbb E_{\mu_\mathcal{G}}\left(\det(\Id_n-G)^{\overline{\delta}}\det(\Id_n-\overline{G})^\delta\right)}.
$$
\end{theorem}
\proof
From Theorem \ref{thm:decompositionRankN}, $G \law R_0\dots R_{n-1}$, hence
\begin{eqnarray*}
\mathbb E_{\mu_{\mathcal{G}}}\left(f(G)\det(\Id_n-G)^{\overline{\delta}}\det(\Id_n-\overline{G})^\delta\right) = \ \ \ \ \ \ \ \ \ \ \ \ \  \ \ \ \ \ \ \ \ \ \ \ \ \ \\ \mathbb E_{\nu_0\star\dots\star\nu_{n-1}}\left(f(R_0\dots
R_{n-1})\det(\Id_n-R_0\dots R_{n-1})^{\overline{\delta}}
\det(\Id_n-\overline{R_0\dots R_{n-1}})^\delta\right)\,.
\end{eqnarray*}
From Lemma \ref{determ}, 
$\det(\Id_n- R_0\dots R_{n-1})=\prod_{k=0}^{n-1}(1 - \langle e_{k+1}, R_k(e_{k+1})\rangle)$, hence
\begin{eqnarray*}
\mathbb E_{\nu_0\star\dots\star\nu_{n-1}}\left(f(R_0\dots
R_{n-1})\det(\Id_n-R_0\dots R_{n-1})^{\overline{\delta}}
\det(\Id_n-\overline{R_0\dots R_{n-1}})^\delta\right)&=&
\\
\mathbb E_{\nu_0\star\dots\star\nu_{n-1}}\left(f(R_0\dots
R_{n-1})\prod_{k=0}^{n-1}\exp_\delta^{(k)}(R_k)\right) \,.\ \ \ \ \ \ \ \ \ \ \ \ \  
\end{eqnarray*}
By the definition of the measures $\nu_k^{(\delta)}$, this is the
desired result.
\qed
\newpage
Before exploring properties of this measure, let us give two 
examples of $\delta$-samplings.

First we check that we can recover
 the classical Ewens measure on the symmetric group. Consider $\mathcal{G}=\mathbb Z_2 \wr {\mathcal S}_n$. For $\delta>0$, the $\delta$-sampling in $\mathbb Z_2 \wr {\mathcal S}_n$ induces a $\theta=2^{2\delta -1}$ sampling on ${\mathcal S}_n$.

\begin{proposition} 
For $\delta > 0$, the pushforward of  $\mu^{(\delta)}_{\mathbb Z_2\wr{\mathcal S_n}}$ by the projection $(f, \sigma) \mapsto \sigma$ is $\mu^{\theta}_n$
with $\theta = 2^{2\delta-1}$.
\end{proposition}

Similarly, if we associate with each transposition of the decomposition (\ref{dec1}) a Rademacher variable, we get easily a sequence of reflections, and if  $\nu_k$ denotes the $k$-th corresponding measure, then the system $(\nu_0, \cdots, \nu_{n-1})$ is coherent with 
$\mu_{\mathbb Z_2\wr{\mathcal S_n}}$. The pushforward of $\nu_k^{(\delta)}$ under the projection is a transposition biased by $\theta$, so we recover the Ewens sampling formula.   

\proof
Recall that the generic element of $\mathbb Z_2 \wr {\mathcal S}_n$ is denoted $(f, \sigma)$. Let ${\mathcal C}(\sigma)$ the set of  cycles of $\sigma$. If $c = (d_1, \dots, d_j)$ is such a cycle, let $\ell(c) =j$ and $w(f; c) = \prod_1^j f(d_j)$. Then it is clear that
\[\det \left(x\Id_n - (f;\sigma)\right) = \prod_{c \in {\mathcal C}(\sigma)} \left(x^{\ell(c)} - w(f; c)\right)\,,\] 
and in particular, 
\begin{equation}
\label{wr1}
\det \left(\Id_n - (f;\sigma)\right)=\left\{
\begin{array}{ccc}
0&\mathrm{if}& \exists c \in {\mathcal C}(\sigma) : w(f;c) = 1\\
2^{k_\sigma}&\mathrm{if}&\forall  c \in {\mathcal C}(\sigma) : w(f;c) = -1 \,.
\end{array}
\right.
\end{equation}
Let $\mathbb P$ stand for $\mu_{\mathbb Z_2\wr{\mathcal S_n}}$ 
 i.e. the uniform distribution on $\mathbb Z_2 \wr {\mathcal S}_n$.
 For any 
test function $F$
\begin{eqnarray*}\mathbb E\left(F(\sigma) |\det (\Id_n - (f, \sigma))|^{2\delta}\right) &=& 
\mathbb E\left[F(\sigma) \mathbb E\left( |\det (\Id_n - (f, \sigma))|^{2\delta} | \sigma\right)\right]\,.
\end{eqnarray*}
Now, conditionally on $\sigma$, the weights of the cycles are independent Rademacher variables (i.e. $\pm 1$ with probability $1/2$). So,
\begin{eqnarray*}
\mathbb P\left( \cap_{c \in {\mathcal C}(\sigma)}\{w(f, \sigma) = -1\}|\sigma\right) = 2^{- k_\sigma}
\end{eqnarray*}
and, due to (\ref{wr1})
\begin{eqnarray*}
\label{wr3}\mathbb E\left( |\det (\Id_n - (f, \sigma))|^{2\delta} | \sigma\right) = 2^{(2\delta- 1) k_\sigma}\,,
\end{eqnarray*}
which  easily yields  
\begin{eqnarray*}
\label{wr2}\mathbb E_{\mu^{(\delta)}_{\mathbb Z_2\wr{\mathcal S_n}}} F(\sigma) = \int_{{\mathcal S}_n} 
F(\sigma) d\mu_n^\theta (\sigma)\,.\end{eqnarray*}
 \qed
 
The fundamental example remains $U(n, \mathbb C)$. In the following section, we will study the determinantal sructure of this model for $\Re\!\ \delta > -1/2$. In \cite{BNR}  a precise analysis of the reflections involved in the decomposition is given. The case $\delta =1$ has a specific interest. If 
 $(\theta_1,\dots,\theta_n)$ are the eigenangles
of a unitary matrix, we have
\[|\det(\Id - U)|^2 = \prod_{j=1}^n |1- \E^{\I\theta_j}|^2\,,\]
which, thanks to the density of the eigenangles, yields 
\begin{multline*}
\mathbb E_{\mu_{U(n)}^{(1)}}\left(f(\theta_1,\dots,\theta_n)\right)\\
= \mbox{cst}\int_{(-\pi,\pi)^n}f(\theta_1,\dots,\theta_n)
\prod_{j<k}|e^{\ii\theta_j}-e^{\ii\theta_k}|^2
\prod_{l=1}^n|1-e^{\ii\theta_l}|^2\D\theta_1\dots\D\theta_n.
\end{multline*}
This means that the distribution of the eigenangles $(\theta_1,\dots,\theta_n)$ of a random matrix drawn according to  $\mu_{U(n)}^{(1)}$ 
 is the same as the distribution of the $n$ first
eigenangles $(\theta_1, \cdots, \theta_n)$ of a random matrix drawn according to $\mu_{U(n+1,\CC)}$, conditionally on $\theta_{n+1}=0$, or,  as seen in \cite{Forr}, 
 as the distribution of $(\theta_1 - \theta_{n+1}, \cdots, \theta_n -\theta_{n+1})$. 
More generally, in \cite{Bourg}, Bourgade gives a geometrical characterisation of this kind of measures for $\delta/2 \in \mathbb N$, defining  
 the notion of conditional Haar measure.

\begin{remark}
A generalized Ewens sampling formula could also be stated
for $\Phi(\mathcal{G})$, with $\mathcal{G}$ checking condition (R)
and $\Phi$ the ring morphism previously defined. 
\end{remark}

\section{A hypergeometric kernel}
\label{Hyper}
In this section, we  study  the correlations of the point process of eigenvalues under the measure $\mu_{U(n,\CC)}^{(\delta)}$ and  answer  Question 3 (see Introduction) asked by Borodin-Olshanski in \cite{BO} section 8. 
Let us recall some basic facts on determinantal processes and correlations, referring
 to the books \cite{AGZ} 4.2 or \cite{Blo} or \cite{ForBook} chap. 4. 
 
Let $\Lambda= \mathbb R$ or $\mathbb T  = \{ z \in \mathbb C : |z|=1\} = \{ \E^{\I \theta} ; \theta \in [-\pi, \pi]\}$ and let us fix an integer $n$.
The collection of eigenvalues $(\lambda_1, \dots, \lambda_n)$ of a random $n\times n$ Hermitian (resp. unitary) matrix can be viewed as  a point process on $\Lambda$ 
,  i.e. a random counting measure $\nu_n = \delta_{\lambda_1} + \dots + \delta_{\lambda_n}$. 
Let us consider  a simple point process $\nu$ on $\Lambda$.
If there exists a sequence of locally integrable functions $\rho_k$ 
 such that for any mutually disjoint family of subsets $D_1, \dots, D_k$ of $\Lambda$ 
\[\mathbb E \left[\prod_{i=1}^k  \nu (D_i)\right] = \int_{\prod_{i=1}^k D_i} \rho_{k} (x_1, \dots, x_k)\D x_1 \dots \D x_k\]
then the functions $\rho_{k}$ are called the correlation functions, or joint intensities of the point process. In this case, the process is said to be determinantal with kernel $K$ if its correlation functions $\rho_k$ are given by 
\[\rho_k(x_1, \dots, x_k) = \det_{i,j=1}^k K(x_i, x_j)\,.\] 
For $\nu= \nu_n$ we denote the correlations by $\rho_{k,n}$ for $k \leq n$. 
When the joint density of the eigenvalues is proportional to
\[\prod_{k=1}^n w(x_k) \prod_{1\leq j<k\leq n} |x_k- x_j|^2\]
for some weight $w$, the orthogonal poynomial method shows that the point process of eigenvalues is determinantal. 
The use of Cayley transform allows to connect Hermitian matrices and unitary matrices.   
 We give a detailed description of the consequence of this connection for the corresponding eigenvalue processes in Subsection \ref{Detproc}, and its impact on the circular Jacobi ensemble in Subsection \ref{our}. Finally, we study the asymptotic behavior in Subsection \ref{limk}.

\subsection{Determinantal processes and Cayley transform}
\label{Detproc}
We follow the approach of Forrester (\cite{ForBook} 2.5 and 4.1.4). 
We start with a weight (positive integrable function) $w^{\mathbb T}$ on $\mathbb T$. The pushforward of the measure
\[\prod_{j=1}w^{\mathbb T}(\E^{\I \theta_j}) \prod_{1\leq j< k\leq n}|\E^{\I\theta_k}- \E^{\I\theta_j}|^2 \D\theta_1 \cdots \D\theta_n\]
by the stereographic projection (Cayley transform)
\[\lambda = \I \frac{1-\E^{\I\theta}}{1+\E^{\I\theta}}= \tan \frac{\theta}{2} \ ; \ \E^{\I\theta}=  \frac{1+\I\lambda}{1-\I\lambda}\]
gives the measure
\[2^{n^2}\prod_{j=1}^n w^{\mathbb T}\left(\frac{1+\I\lambda_j}{1-\I\lambda_j}\right) (1+\lambda_j^2)^{-n} \prod_{1\leq j< k\leq n}|\lambda_k- \lambda_j|^2 \D\lambda_1 \cdots \D\lambda_n.\]
We define the weight $w^{\mathbb R}$ on $\mathbb R$ as
\[w^{\mathbb R} (x) = (1+x^2)^{-n}w^{\mathbb T}\left(\frac{1+\I x}{1-\I x}\right)\,.\]
Conversely
\[w^{\mathbb T}(\E^{\I\theta}) = \left(\cos \frac{\theta}{2}\right)^{2n} w^{\mathbb R}\left(\tan \frac{\theta}{2}\right)\,.\]
If the monomials $1, x, \dots, x^n$ are in $L^2(w^{\mathbb R} (x)\D x)$, then the orthogonal polynomial method gives
\[\frac{1}{{\mathcal Z}_n^{\mathbb R}}\prod_{j=1}^n w^{\mathbb R}(\lambda_j)
\prod_{1\leq j< k\leq n} |\lambda_k- \lambda_j|^2 = \det \big(\widetilde K^{\mathbb R}_n (\lambda_j, \lambda_k)\big)_{1 \leq j,k\leq n} \]
where $\mathcal Z_n^{\mathbb R}$ is a normalization constant and where
\begin{eqnarray*}\widetilde K^{\mathbb R}_n (x,y) &=& \sqrt{w^{\mathbb R} (x)w^{\mathbb R} (y)}\ K^{\mathbb R}_n (x,y)\\ 
K^{\mathbb R}_n (x,y)&=& 
\sum_{\ell =0}^{n-1} p_\ell^{\mathbb R}(x)p_\ell^{\mathbb R}(y) \end{eqnarray*}
and the $p^{\mathbb R}_\ell$ are orthonormal with respect to 
 the measure $w^{\mathbb R}(x) \D x$.
The Christoffel-Darboux formula gives another expression for the kernel
\[K^{\mathbb R}_n (x,y) = \frac{\kappa_{n-1}}{\kappa_n}\ \frac{p^{\mathbb R}_{n}(x)p^{\mathbb R}_{n-1}(y) - p^{\mathbb R}_{n-1}(x)p^{\mathbb R}_{n}(y)}{x-y}\] 
where $\kappa_j$ is the coefficient of $x^j$ in $p_j^{\mathbb R}(x)$.
In terms of the monic orthogonal polynomials $P_0, \cdots, P_{n}$, this yields
\begin{eqnarray}
K_n^{\mathbb R} (x,y) &=& \sum_{\ell = 0}^{n-1} \frac{P_\ell (x)P_\ell (y)}{\Vert P_\ell\Vert^2}\\
&=& 
 \frac{P^{\mathbb R}_{n}(x)P^{\mathbb R}_{n-1}(y) - P^{\mathbb R}_{n-1}(x)P^{\mathbb R}_{n}(y)}{\Vert P_{n-1}\Vert^2 (x-y)}\,.
\end{eqnarray}

Besides, on the unit circle, we consider the polynomials $\varphi_\ell$ (resp. $\Phi_\ell$) orthonormal (resp. monic orthogonal) with respect to  the measure $w^{\mathbb T}(\E^{\I\theta})\D\theta$, and their reciprocal defined by  
\[ \Phi_\ell^\star(z) = z^\ell\!\ \overline{\Phi_\ell (1/ \bar z)}\ , \ \varphi_\ell^\star(z) = z^\ell\!\ \overline{\varphi_\ell (1/ \bar z)}\,.\]
We have then 
\[\frac{1}{{\mathcal Z}_n^{\mathbb T}} \prod_{j=1}^n w^{\mathbb T}(\E^{\I\theta_j}) \prod_{1\leq j< k\leq n}|\E^{\I\theta_k}- \E^{\I\theta_j}|^2 = \det \big(\widetilde K^{\mathbb T}_n (\E^{\I\theta_j}, \E^{\I\theta_k})\big)_{1 \leq j,k\leq n} \]
with
\[\widetilde K^{\mathbb T}_n (z,\zeta) = \sqrt{w^{\mathbb T} (z)w^{\mathbb T} (\zeta)}\ 
 K^{\mathbb T}_n (z,\zeta) \]
 and
 \[K^{\mathbb T}_n (z,\zeta)= 
\sum_{\ell =0}^{n-1} \overline{\varphi_\ell(z)}\varphi_\ell (\zeta)\,.\]
%
\newpage
The Christoffel-Darboux formula is now
\begin{equation}
\label{CDC}K_n^{\mathbb T}(z, \zeta)= \frac{\overline{\varphi_{n}^*(z)}\varphi_{n}^*(\zeta) - \overline{\varphi_{n}(z)}\varphi_{n}(\zeta)}{1-\bar z\zeta}\end{equation}
(see \cite{SimonCD} 1.12 and 3.2), or 
\begin{equation}
K_n^{\mathbb T}(z, \zeta) = 
\frac{\overline{\Phi_{n}^*(z)}\Phi_{n}^*(\zeta) - \overline{\Phi_{n}(z)}\Phi_{n}(\zeta)}{\Vert \Phi_n \Vert^2(1-\bar z\zeta)}\,.\end{equation}
The kernel $\widetilde K_n^{\mathbb R}$ (resp. $\widetilde K_n^{\mathbb T}$) rules the correlation function 
$\rho_{n,m}^{\mathbb R}(\lambda_1, \cdots, \lambda_m)$ (resp. $\rho_{n,m}^{\mathbb C}(\E^{\I\theta_1}, \cdots, \E^{\I\theta_m})$) for $m=1, \cdots, n$.

\subsection{Our weights and their characteristics}
\label{our}
For the sake of simplicity we use the polygamma symbol
\[\Gamma\left[\begin{matrix}\let \over/ a, b, \cdots\\ \let \over /c,d, \cdots\end{matrix}\right] := \frac{\Gamma(a)\Gamma(b)\cdots}{\Gamma(c)\Gamma(d)\cdots}.\]

For $\delta = a+\I b\in \mathbb C$ with $a > -1/2$, we will consider two weights on $(-\pi, \pi)$
\begin{eqnarray}
\label{2w}
w_1^{\mathbb T}(\E^{\I\theta}) &=& (1-\E^{\I\theta})^{\overline{\delta}}(1-\E^{-\I\theta})^{\delta} = (2-2\cos\theta)^a \E^{-b(\pi \sgn\theta-\theta)}\\
w_2^{\mathbb T}(\E^{\I\theta}) &=& (1+\E^{\I\theta})^{\overline{\delta}}(1+\E^{-\I\theta})^{\delta} = (2+2\cos\theta)^a \E^{-b\theta}
\end{eqnarray}
These are "pure" Fisher-Hartwig functions.
We can go from $w_1^{\mathbb T}$ to $w_2^{\mathbb T}$  by the transform
\begin{equation}\label{thetatau}\theta \mapsto \tau:= -\theta +\pi(\sgn \theta)\end{equation}
which carries the discontinuity in $\theta=0$ to the edges $\pm \pi$, so that
\begin{equation}
\label{thetatau1}\E^{\I\theta} = - \E^{-\I\tau} \ \hbox{and} \ w_1^{\mathbb T}(\E^{\I\theta}) = w_2^{\mathbb T}(\E^{-\I\tau})\,.\end{equation}

For $a > -1/2$, the Fourier coefficients of  $w_1$ are known (\cite{BoS} Lemma 2.1)
\[\frac{1}{2\pi}\int_{-\pi}^{\pi} w_1^{\mathbb T} (\E^{\I\theta}) \E^{-\I n\theta} \D\theta = (-1)^n 
\Gamma\left[\begin{matrix}\let \over/ 1+\delta+\bar\delta\\  \let \over / \bar\delta -n +1 , \delta+n+1 \end{matrix}\right].\]
With$$ c(\delta)= \frac{1}{2\pi}\Gamma\left[\begin{matrix}\let \over/ 1+\delta, 1+\bar\delta \\ \let \over /1 + \delta  + \bar\delta\end{matrix}\right],$$
 the  function $\widetilde w_1^{\mathbb T} (\E^{\I\theta}) = c(\delta)w_1^{\mathbb T}(\E^{\I\theta})$ is a  probability density on $(-\pi, \pi)$.
For $w_2$, we note that
\[\int_{-\pi}^\pi w_1^{\mathbb T}(\E^{\I\theta})  \E^{-\I n\theta} \D\theta = (-1)^n \int_{-\pi}^\pi w_2^{\mathbb T}(\E^{\I\tau})  \E^{\I n\tau} \D\tau.\]
Moreover we go from one system of polynomials to the other by the mapping $z \mapsto -z$.
\newpage

It is known from \cite{Askey} p. 304 and \cite{Basor} p.31-34 that for $n \geq 0$ the  $n$-th orthonormal polynomial with respect to $\widetilde w_1^{\mathbb T}(\E^{\I\theta}) \D\theta$ is
\begin{eqnarray}
\label{vraiphi}\Phi_n (z) = 
\Gamma\left[\begin{matrix}\let \over/\delta + n, \bar\delta + 1\\ \let \over /\bar\delta + n+1, \delta\end{matrix}\ \right]
\ _2F_1\left(\begin{matrix}\let \over/-n, \bar\delta+1\\ \let \over /1-n-\delta\end{matrix}\ ;\ z\right)
\end{eqnarray}
with
\begin{equation}
\label{normphi}\Vert \Phi_n\Vert^2 = 
\Gamma\left[\begin{matrix}\let \over/\delta +\bar\delta + n+1 , n+1, \bar\delta + 1, \delta +1 \\ \let \over /\bar\delta + n+1, \delta + n+1, \delta + \bar\delta +1\end{matrix}\ \right]\,,\end{equation}
(see also \cite{ForBook} Prop. 4.8 in the case $\delta$ real).
With the complement formula (\ref{1-z}) we get the other form
\begin{eqnarray}
\label{vraiphicool}
\Phi_n (z) = \Gamma\left[\begin{matrix}\let \over/\delta+ \bar\delta +1+n, \bar\delta + 1\\ \let \over /\bar\delta +n +1, \delta + \bar\delta +1 \end{matrix}\ \right]\   _2F_1\left(\begin{matrix}\let \over/-n, \bar\delta + 1\\ \let \over /\delta + \bar\delta +1\end{matrix}\ ;\ 1-z\right).
\end{eqnarray}
 In view  of (\ref{1/z}) and (\ref{vraiphi}) we identify $\Phi^*_n$ 
as
\begin{eqnarray}
\label{vraiphistar}
\Phi^*_n (z) = \ _2F_1\left(\begin{matrix}\let \over/-n, \bar\delta\\ \let \over /-n-\delta\end{matrix}\ ;\ z\right),\,
\end{eqnarray}
or, using (\ref{1-z}) again
\begin{equation}
\label{vraiphistarcool}
\Phi^*_n (z) = \Gamma\left[\begin{matrix}\let \over/\delta+ \bar\delta +1+n, \delta + 1\\ \let \over /\delta +n +1, \delta + \bar\delta +1 \end{matrix}\ \right]\ _2F_1\left(\begin{matrix}\let \over/-n, \bar\delta\\ \let \over /\delta + \bar\delta +1\end{matrix}\ ;\ 1 -z\right).
\end{equation}
Borodin and Olshanski  considered the following weight on $\mathbb R$ : 
\begin{equation}
\label{BOw2}2^{-\delta -\bar\delta} w_2^{\mathbb R}(x)
 = (1+\I x)^{-\delta-n}(1-\I x)^{-\bar\delta -n}\,.\end{equation}
Since this weight depends on $n$, the reference measure has only a finite set of moments
so that there is only a finite set of orthogonal polynomials  (these are the pseudo-Jacobi polynomials)
\begin{equation}p_m (x) = (x-\I)^m \ _2F_1\left(\begin{matrix}\let \over/-m,\delta +n-m\\ \let \over /\delta+\bar\delta +2n-2m\end{matrix}\ ;\ \frac{2}{1+\I x}\right)\end{equation}
$m < a+n - \frac{1}{2}$.
Let us call $\widetilde K_{2,n}^{\mathbb R}$ the corresponding kernel.
\subsection{Asymptotic behavior}
\label{limk}
For the weight $w_2^{\mathbb R}$, Borodin and Olshanski  considered the (thermodynamic) scaling limit $\lambda \mapsto n\lambda$ and proved  (\cite{BO} Theorem 2.1)
\begin{theorem}[Borodin-Olshanski]
Let $\Re\!\ \delta > -1/2$. 
\begin{enumerate}
\item We have
\begin{eqnarray}
\label{BOlim}\lim_n\ (\sgn\!\ x\!\ \sgn\!\ y)^n n\widetilde K_{2,n}^{\mathbb R} (nx,ny) &=& \widetilde K_\infty^{\mathbb R} (x,y)
\end{eqnarray}
uniformly for $x, y$ in compact sets of $\mathbb R^\star \times \mathbb R^\star$,
where (for $x \not= y$)
\begin{eqnarray}
\label{defkr}
\widetilde K_\infty^{\mathbb R} (x,y)
&:=& \frac{1}{2\pi}\Gamma\left[\begin{matrix}\let \over/\delta +1, \bar\delta+1\\ \let \over /\delta + \bar\delta +1, \delta + \bar\delta +2\end{matrix}\right] \frac{\widetilde P(x)Q(y)- Q(x)\widetilde P(y)}{x-y} \\
\label{wherep}
\widetilde P(x) &=& \left|\frac{2}{x}\right|^{\frac{\delta + \bar\delta}{2}}\E^{-\frac{\I}{x}+\pi\frac{(\delta-\bar\delta)\sgn x}{4}}\ _1F_1\left(\begin{matrix}\let \over/\delta\\ \let \over /\delta+ \bar\delta +1\end{matrix}\ ;\ \frac{2\I}{x}\right)\\ 
 \label{whereq} 
Q(x) &=& \frac{2}{x}\left|\frac{2}{x}\right|^{\frac{\delta + \bar\delta}{2}}\E^{-\frac{\I}{x}+\pi\frac{(\delta-\bar\delta)\sgn x}{4}}\ _1F_1\left(\begin{matrix}\let \over/\delta +1\\ \let \over /\delta+ \bar\delta+2\end{matrix}\ ;\ \frac{2\I}{x}\right)\,. 
\end{eqnarray}
\item The limiting correlation is given by
\begin{equation}
\label{limcorrR}\lim_n n^m \rho^{\mathbb R}_{n,m}(n\lambda_1, \cdots, n\lambda_m) = \det\left(\widetilde K_\infty^{\mathbb R}(\lambda_i, \lambda_j)\right)_{1\leq i, j\leq m}\,.\end{equation}
\end{enumerate}
\end{theorem}
The kernel $\widetilde K_\infty^{\mathbb R}(1/x, 1/y)$ is called the \textsl{confluent hypergeometric kernel} in \cite{BD}.

For the circular model, we choose the set-up $w_1$ for the sake of consistency with the above sections. The singularity is in $z=1$ i.e. $\theta= 0$. To study the asymptotic behavior of the point process on $\mathbb T$ at the singularity (edge) we have two ways: either
 take the thermodynamic scaling $\theta \mapsto \theta/n$, or 
 use the result on $\mathbb R$.

\begin{theorem}
\label{theokernel}
Let $\Re\!\ \delta > -1/2$. 
\begin{enumerate}
\item
With the weight $w_1$,  
\begin{eqnarray}
\label{cvK}
\lim_n n^{-1} \widetilde K^{\mathbb T, 1}_n (\E^{\I\theta/n}, \E^{\I\tau/n}) = \widetilde K^{\mathbb T}_\infty(\theta,\tau)
\end{eqnarray}
with, for $\theta\not=\tau$
\begin{eqnarray}
\nonumber
\widetilde K^{\mathbb T}_\infty(\theta,\tau)
&=& 
\frac{1}{2\I\pi}
\Gamma\left[\begin{matrix}\let \over/1+\delta, 1+\bar\delta\\ \let \over /1+\delta+\bar\delta, 1+\delta+\bar\delta\end{matrix}\ \right]\frac{P^{\mathbb T}(\theta)\overline{P^{\mathbb T}(\tau)}- \overline{P^{\mathbb T}(\theta)}P^{\mathbb T}(\tau)} {\theta-\tau}\\
\label{defkc}
\end{eqnarray}
where
\begin{eqnarray}
\nonumber 
P^{\mathbb T}(\theta) := |\theta|^{\frac{\delta + \bar\delta}{2}}\E^{\I \frac{\theta}{2}- \frac{\pi}{4}(\delta-\bar\delta)\sgn \theta}\ 
_1F_1\left(\begin{matrix}\let \over/\delta\\ \let \over /\delta + \bar\delta +1\end{matrix}\ ;\ \displaystyle{-\I \theta}\right) = \widetilde P\big( -2\theta^{-1}\big)\,,\\
\label{PPtilde}
\end{eqnarray}
 and
\begin{eqnarray}
\nonumber
\widetilde K^{\mathbb T}_\infty(\theta,\theta) = \frac{|\theta|^{\delta + \bar\delta}}{2\pi}\Gamma\left[\begin{matrix}\let \over/1+\delta, 1+\bar\delta\\ \let \over /1+\delta+\bar\delta, 1+\delta+\bar\delta\end{matrix}\ \right]
 \Re\ \ \ \ \ \ \ \ \ \ \ \ \ \ \ \ \ \ \ \ \ \ \ \ \ \ \ \ \ \ \ \ \ \ \ \ \ \ \ \ \ \ \ \ \ \ \ \ \ \ \ \ \ \ \ \ 
\end{eqnarray}
\begin{eqnarray}
 \nonumber
\left[  _1F_1\left(\begin{matrix}\let \over/\delta\\ \let \over /\delta + \bar\delta +1\end{matrix}\ ; -\I \theta\right)
\left[_1F_1\left(\begin{matrix}\let \over/\bar\delta\\ \let \over /\delta + \bar\delta +1\end{matrix}\ ; \I \theta\right) -2\  _1F_1\left(\begin{matrix}\let \over/\bar\delta +1\\ \let \over /\delta + \bar\delta +2\end{matrix}\ ; \I \theta\right)
\right]\right]
\\
\end{eqnarray}
\item 
The limiting correlation is given by
\begin{equation}
\label{limcorrT}\lim_n n^m \rho^{\mathbb T, 1}_{n,m}(\E^{\I\theta_1/n}, \cdots, \E^{\I\theta_m/n}) = \det\left(\widetilde K_\infty^{\mathbb T}(\theta_i, \theta_j)\right)_{1\leq i, j\leq m}\,.\end{equation}
\end{enumerate}
\end{theorem}

\proof We begin with a direct proof of (\ref{cvK}) when $\theta\not=\tau$, and then  proceed with the proof of (\ref{cvK}) when $\theta=\tau$, which directly yields  (\ref{limcorrT}) and  we end  with an alternate proof of (\ref{limcorrT}) using (\ref{limcorrR}) and the Cayley transform. 
\medskip

\noindent 1) 
The following lemma describes  the asymptotical behavior of the quantities entering in the kernel.
\begin{lemma}
When $n\rightarrow \infty$ 
\begin{equation}
\label{normlim}
\lim_n \Vert \Phi_n\Vert^2 = \Gamma\left[\begin{matrix}\let \over/\bar\delta +1 , \delta + 1\\ \let \over /\delta + \bar\delta +1\end{matrix}\ \right] 
\end{equation} 
Moreover if  $n\theta_n \rightarrow \theta$, then (uniformly for $\theta$ in a compact set)
\begin{eqnarray}\label{asphi}
\lim n^{-\delta}\Phi_n ( \E^{\I\theta_n}) &=& \Gamma\left[\begin{matrix}\let \over/\bar\delta +1 \\ \let \over / \delta + \bar\delta +1\end{matrix}\ \right] 
\ _1F_1\left(\begin{matrix}\let \over/\bar\delta+1\\ \let \over /\delta + \bar\delta +1\end{matrix}\ ;\ \displaystyle{\I \theta}\right),
 \\ \label{asphistar}
\lim n^{-\bar\delta}\Phi^\star_n ( \E^{\I\theta_n}) &=& \Gamma\left[\begin{matrix}\let \over/\delta +1 \\ \let \over / \delta + \bar\delta +1\end{matrix}\ \right] 
\ _1F_1\left(\begin{matrix}\let \over/\bar\delta\\ \let \over /\delta + \bar\delta +1\end{matrix}\ ;\ \displaystyle{\I \theta}\right),
\\ \label{asphi'}
\lim n^{-\bar\delta +1}(\Phi^\star_n)' ( \E^{\I\theta_n}) &=& \bar\delta\Gamma\left[\begin{matrix}\let \over/\delta +1 \\ \let \over / \delta + \bar\delta +2\end{matrix}\ \right]\ _1F_1\left(\begin{matrix}\let \over/ \bar\delta+1\\ \let \over /\delta+ \bar\delta + 2\end{matrix}\ ;\ \I\theta\right).
\end{eqnarray}
\end{lemma}

\proof 
Let us first recall that, as $n \rightarrow \infty$,
\begin{equation}\label{pot}\frac{\Gamma(c+n)}{\Gamma(n)}\sim n^{c}\,,\end{equation}
which gives immediately (\ref{normlim}).
The limits in (\ref{asphi}) and (\ref{asphistar}) are then consequences of (\ref{vraiphicool}), (\ref{vraiphistarcool}) and the limiting relation (\ref{cvcfl}).
Besides, in view of (\ref{deriv}) and (\ref{vraiphistarcool}),
\[(\Phi^\star_n)' (z) = \frac{n\bar\delta}{\delta+\bar\delta +1} \Gamma\left[\begin{matrix}\let \over/\delta+\bar\delta +n , \delta +1 \\ \let \over / \delta + \bar\delta +1 , \delta +n +1\end{matrix}\ \right] \ _2F_1\left(\begin{matrix}\let \over/-n+1, \bar\delta+1\\ \let \over /\delta+\bar\delta + 2\end{matrix}\ ;\ 1-z\right)\,.\]
It remains to apply 
(\ref{cvcfl}). \qed
\medskip

\noindent\underline{A) For $\theta \not= \tau$}, we have, by the Christoffel-Darboux formula (\ref{CDC}):
\begin{eqnarray*}\lim_n \!\ \I (\theta-\tau) \Gamma(\delta + \bar\delta +1) n^{-(\delta + \bar \delta +1)} K_n^{\mathbb T, 1}(\E^{\I\theta/n}, \E^{\I\tau/n}) =\\
_1F_1\left(\begin{matrix}\let \over/\delta\\ \let \over /\delta + \bar\delta +1\end{matrix}\ ;\ \displaystyle{-\I \theta}\right)\!\  _1F_1\left(\begin{matrix}\let \over/\bar\delta\\ \let \over /\delta + \bar\delta +1\end{matrix}\ ;\ \displaystyle{\I \tau}\right) \\-\!\  _1F_1\left(\begin{matrix}\let \over/\delta+1\\ \let \over /\delta + \bar\delta +1\end{matrix}\ ;\ \displaystyle{-\I \theta}\right)\ 
_1F_1\left(\begin{matrix}\let \over/\bar\delta +1\\ \let \over /\delta + \bar\delta +1\end{matrix}\ ;\ \displaystyle{\I \tau}\right)\end{eqnarray*}
Now, applying  the Kummer's formula (\ref{1F1m})
\begin{eqnarray*}_1F_1\left(\begin{matrix}\let \over/\delta+1\\ \let \over /\delta + \bar\delta +1\end{matrix}\ ;\ \displaystyle{-\I \theta}\right)
&=&\E^{-\I\theta}\ _1F_1\left(\begin{matrix}\let \over/\bar\delta\\ \let \over /\delta + \bar\delta +1\end{matrix}\ ;\ \displaystyle{\I \theta}\right)\\
\ 
_1F_1\left(\begin{matrix}\let \over/\bar\delta +1\\ \let \over /\delta + \bar\delta +1\end{matrix}\ ;\ \displaystyle{\I \tau}\right)
&=& \E^{\I\tau}\ 
_1F_1\left(\begin{matrix}\let \over/\delta\\ \let \over /\delta + \bar\delta +1\end{matrix}\ ;\ \displaystyle{-\I \tau}\right)
\end{eqnarray*}
Besides we have (recall that we used $\widetilde w_1$)
\[\frac{\widetilde K_n^{\mathbb T, 1}(\E^{\I\theta/n}, \E^{\I\tau/n}) }{K_n^{\mathbb T, 1}(\E^{\I\theta/n}, \E^{\I\tau/n}) }
= c(\delta)\sqrt{w_1(\E^{\I\theta/n})w_1(\E^{\I\tau/n})}\]
and from the very definition of $w_1$
\[\lim n^{2(\delta + \bar\delta)}w_1(\E^{\I\theta/n})w_1(\E^{\I\tau/n}) = |\theta\tau|^{2\Re \delta} \E^{-\Im\delta \pi(\sgn \theta + \sgn \tau)}\]
We conclude that (\ref{cvK}) holds true.

\medskip
\noindent\underbar{B) On the diagonal} In the following $z$ and $\zeta$ are elements of $\mathbb T$. If $F$ and $G$ are differentiable functions on $\mathbb T$, the  de l'Hospital rule gives
\[\lim_{\zeta \rightarrow z} \frac{F(z) G(\zeta) -F(\zeta) G(z)}{z-\zeta} = F'(z)G(z) -F(z)G'(z)\,.\]
Taking 
\[F(z) = z^{-n}\Phi_n(z) \ , \ G(z) = \overline{\Phi_n(z)}\,,\]
so that
\[ F'(z) = -n z^{-n-1}\Phi_n(z) + z^{-n}\Phi_n'(z)\ , \ G'(z) = -z^{-2} \overline{\Phi_n'(z)}\]
we get the value of the kernel on the diagonal:
\begin{eqnarray}\nonumber\lim_{\zeta \rightarrow z} \frac{\overline{\Phi_n^*(z)}\Phi_n^*(\zeta) - \overline{\Phi_n(z)}\Phi_n(\zeta)}{1 - \bar z \zeta} &=& - n|\Phi_n(z)|^2 + 2 \Re[\overline{\Phi_n(z)} z\Phi_n'(z)]\\
\nonumber
&=& n|\Phi^*_n(z)|^2 - 2 \Re[\overline{\Phi^*_n(z)} z(\Phi^*_n)'(z)].\\
\end{eqnarray}
It remains to apply the lemma. 

Notice that
\[\lim_n n^{-(1+\delta + \bar\delta)} K_n^{\mathbb T, 1} (1, 1) = \frac{1}{\Gamma(\delta + \bar \delta + 2)}\,.\]

\medskip
\noindent 2) \textsl{Alternate proof of (\ref{limcorrT})}

The pushforward of the measure 
\[\rho_{n}^{\mathbb R, 2}(x_1, \cdots, x_n) dx_1 \dots dx_n\]
by  the Cayley transform is, 
\[2^{-n}\rho_{n}^{\mathbb R, 2}\left(\tan\frac{\theta_1}{2},\cdots,\tan\frac{\theta_n}{2}\right)
\prod_{k=1}^n \cos^{-2} \frac{\theta_k}{2}\ \D\theta_1 \dots \D\theta_n\]
which, at the level of kernels gives
\[\rho_{n, m}^{\mathbb T, 2}(\E^{\I\theta_1},\cdots, \E^{\I\theta_m}) ) = \det\left[\widetilde K_n^{\mathbb R, 2}\left(\tan \frac{\theta_i}{2}, \tan \frac{\theta_j}{2}\right) \frac{1}{2 \cos\theta_i\cos\theta_j}\right]_{1\leq i, j \leq m}. \]
Coming back to the superscript $1$ with the help of (\ref{thetatau}) we obtain
\[\rho_{n, m}^{\mathbb T, 1}(\E^{\I\theta_1},\cdots, \E^{\I\theta_m}) ) = \det\left[ H_n(\theta_i), H_n (\theta_j)\right]_{1\leq i, j \leq m} \]
with
\[H_n(\theta, \theta') = \widetilde K _n ^{\mathbb R, 2}\left(- \cot \frac{\theta}{2}, -\cot\frac{\theta'}{2}\right) \frac{1}{2 |\sin\frac{\theta}{2}\sin\frac{\theta'}{2}|}.\]

Let us rescale the angles. Since $\lim_n n \tan\frac{\theta}{n} = \theta$ , $\lim_n n \tan\frac{\theta'}{n} = \theta'$ and since the limit in (\ref{BOlim}) is uniform on compact subsets, we get
\[\lim \frac{1}{n} H_n \left(\frac{\theta}{n} , \frac{\theta'}{n}\right) = \frac{2}{|\theta\theta'|} \widetilde K_\infty^{\mathbb R} \left(-\frac{2}{\theta}, -\frac{2}{\theta'}\right)\,.\]

We remark that  $P^{\mathbb T}(\theta)=\widetilde P(x)$ with $x\theta= -2$. Moreover, from (\ref{recur}), we have
\[\frac{\I}{\bar\delta + \delta +1} Q(x) = \overline{P^{\mathbb T}(\theta)} - P^{\mathbb T}(\theta)\]
so that, if $\tau = -2/y$
\[\frac{\I}{\bar\delta + \delta +1}\left[\widetilde P(x) Q(y) - \widetilde P(y) Q(x)\right] = P^{\mathbb T}(\theta)\overline{P^{\mathbb T}(\tau)} - 
P^{\mathbb T}(\tau)\overline{P^{\mathbb T}(\theta)}\] 
and consequently
\begin{equation}
\frac{\theta\tau}{2} \widetilde K_\infty^{\mathbb T}(\theta, \tau) = \widetilde K_\infty^{\mathbb R}(x, y)\,.
\end{equation}
\qed
\begin{remark}
\begin{enumerate}
\item
To  have a graphical point of view of this kernel, we refer to \cite{BourgT} p.56--60.
\item In \cite{NNR},  the behavior of the limiting kernel on $\mathbb R$ is used to study asymptotics of the maximal eigenvalue of the generalized Cauchy ensemble.
\item An easy computation shows that for $\delta$ real, $\delta>-1/2$, we recover the Bessel kernel
$$K_\infty^{\mathbb{T}}=\dfrac{\pi}{2}\sqrt{\theta \tau} \dfrac{J_{\delta+\frac{1}{2}}(\frac{\pi \theta}{2})J_{\delta-\frac{1}{2}}(\frac{\pi \tau}{2})-J_{\delta-\frac{1}{2}}(\frac{\pi \theta}{2})J_{\delta+\frac{1}{2}}(\frac{\pi \tau}{2})}{2(\theta-\tau)},$$
and for $\delta=0$ the sine kernel
$$K_\infty^{\mathbb{T}}=\dfrac{\sin(\frac{\theta-\tau}{2})}{\pi(\theta-\tau)}.$$
\end{enumerate}
\end{remark}

\section{Appendix: Hypergeometric functions}
For a classical reference on hypergeometric functions, see \cite{AAR}. 

The Gauss hypergeometric function is defined as
\begin{equation}
\label{defhypf}
_2F_1\left(\begin{matrix}\let \over/a,b\\ \let \over /c\end{matrix}\ ;\ z\right) = \sum_{k=0}^\infty \frac{(a)_k (b)_k}{(c)_k} \frac{z^k}{k!}
\end{equation}
where $(x)_n$ stands for the
Pochhammer symbol  $(x)_k =x(x+1)\dots(x+k-1)$, with the convention $(x)_0 = 1$. 
When $a=-n\in -\mathbb N_0$, it is a polynomial 
\begin{equation}
\label{defhypp}
_2F_1\left(\begin{matrix}\let \over/-n,b\\ \let \over /c\end{matrix}\ ;\ z\right) = \sum_{k=0}^n (-1)^k \binom{n}{k}\frac{(b)_k}{(c)_k} z^k.
\end{equation}
The following relations are useful:
\begin{equation}
\label{1/z}z^n\ _2F_1\left(\begin{matrix}\let \over/-n,b\\ \let \over /c\end{matrix}\ ;\ z^{-1}\right) = (-1)^n \frac{(b)_n}{(c)_n}\ _2F_1
\left(\begin{matrix}\let \over/-n,-n-c+1\\ \let \over /-n-b+1\end{matrix}\ ;\ z\right)
\end{equation}
\begin{equation}
\label{1-z}
_2F_1\left(\begin{matrix}\let \over/-n,b\\ \let \over /c\end{matrix}\ ;\ \displaystyle{1-z}\right) = \frac{(c-b)_n}{(c)_n}\ _2F_1\left(\begin{matrix}\let \over/-n,b\\
\let\over/ -n+b+1-c\end{matrix} \ ;\ \displaystyle{z}\right)\end{equation}
\begin{equation}
\label{deriv}
\frac{d}{dz}\ _2F_1\left(\begin{matrix}\let \over/a,b\\ \let \over /c\end{matrix}\ ;\ \displaystyle z\right) =  \frac{ab}{c} \ _2F_1\left(\begin{matrix}\let \over/a+1,b+1\\ \let \over /c+1\end{matrix}\ ;\ \displaystyle z\right)
\,.\end{equation}
It is known that, uniformly for $z$ in a compact set, for $b,c$ fixed 
\begin{equation}
\label{cvcfl}
\lim_N \ _2F_1\left(\begin{matrix}\let \over/-N,b\\ \let \over /c\end{matrix}\ ;\ -\frac{z}{N}\right) =\ _1F_1\left(\begin{matrix}\let \over/b\\ \let \over /c\end{matrix}\ ;\ \displaystyle{z}\right)
\end{equation}
where 
\begin{equation}
\label{confluent}_1F_1 \left(\begin{matrix}\let \over/b\\ \let \over /c\end{matrix}\ ;\ \displaystyle{z}\right) =\sum_{k=0}^\infty \frac{(b)_k}{(c)_k} \frac{z^k}{k!}\end{equation}
is the confluent hypergeometric function.

It satisfies  Kummer's formula: 
\begin{equation}
\label{1F1m}
\E^z \ _1F_1\left(\begin{matrix}\let \over/a\\ \let \over /c\end{matrix}\ ;\ \displaystyle{-z}\right) = \ _1F_1\left(\begin{matrix}\let \over/c-a\\ \let \over /c\end{matrix}\ ;\ \displaystyle{z}\right),\end{equation}
the recursion formula
\begin{equation}
\label{recur}
_1F_1\left(\begin{matrix}\let \over/a\\ \let \over /c\end{matrix}\ ;\ \displaystyle{z}\right)=\ _1F_1\left(\begin{matrix}\let \over/a-1 \\ \let \over /c\end{matrix}\ ;\ \displaystyle{z}\right)+ \frac{z}{c}\ _1F_1\left(\begin{matrix}\let \over/a\\ \let \over /c+1\end{matrix}\ ;\ \displaystyle{z}\right)\,,  
\end{equation}
and the derivative formula
\begin{equation}
\label{derivconf}
\frac{d}{dz}\ _1F_1\left(\begin{matrix}\let \over/a\\ \let \over /c\end{matrix}\ ;\ \displaystyle z\right) =  \frac{a}{c} \ _1F_1\left(\begin{matrix}\let \over/a+1\\ \let \over /c+1\end{matrix}\ ;\ \displaystyle z\right)
\,.
\end{equation}
\bigskip

\noindent{\bf Acknowledgement}
A.N.'s work is supported by the Swiss National Science Foundation (SNF) grant 200021\_119970/1.

A.R's work is partly supported by the ANR project Grandes Matrices
Al\'eatoires ANR-08-BLAN-0311-01.

\renewcommand{\refname}{References}

\end{document}